\documentclass{amsart}

\usepackage{amsmath,amsfonts,amssymb,amscd}

\usepackage{isolatin1}

\usepackage[matrix,frame]{xypic}

\newtheorem{e-proposition}[theorem]{Proposition}

\newtheorem{e-definition}[theorem]{Definition\rm}

\newtheorem{theoreme}{Th\'eor\`eme}[section]
\newtheorem{lemme}[theoreme]{Lemme}
\newtheorem{proposition}[theoreme]{Proposition}
\newtheorem{corollaire}[theoreme]{Corollaire}

\newtheorem{conjecture}[theoreme]{Conjecture}

\def\og{\leavevmode\raise.3ex\hbox{$\scriptscriptstyle\langle\!\langle$~}}
\def\fg{\leavevmode\raise.3ex\hbox{~$\!\scriptscriptstyle\,\rangle\!\rangle$}}

\def\<{\langle\,}
\def\>{\,\rangle}

\def\shuff#1#2{\mathbin{
      \hbox{\vbox{
        \hbox{\vrule
              \hskip#2
              \vrule height#1 width 0pt
               }%
        \hrule}%
             \vbox{
        \hbox{\vrule
              \hskip#2
              \vrule height#1 width 0pt
               \vrule }%
        \hrule}%
}}}
\def\shuffl{{\mathchoice{\shuff{7pt}{3.5pt}}%
                        {\shuff{6pt}{3pt}}%
                        {\shuff{4pt}{2pt}}%
                        {\shuff{3pt}{1.5pt}}}}%
\def\shuffle{\, \shuffl \,}

\def\SG{{\mathfrak S}}

\def\FQSym{{\bf FQSym}}

\def\F{{\bf F}}

\def\Z{{\mathbb Z}}

\newcommand{\sign}{{\operatorname{sgn}}}
\newcommand{\del}{{\operatorname{del}}}

\newcommand{\Des}{{\operatorname{Des}}}
\newcommand{\Rec}{{\operatorname{Rec}}}
\newcommand{\CDes}{{\operatorname{C}}}

\begin{document}




%
\title{Sur une conjecture de Dehornoy}

\vspace{-2.6cm}
\title{On a conjecture by Dehornoy}



\author[F. Hivert, J.-C.~Novelli, and J.-Y.~Thibon]%
{Florent Hivert, Jean-Christophe Novelli, and Jean-Yves Thibon}

\address[Hivert]{LITIS, Universit\'e de Rouen ; Avenue de l'universit\'e ;
76801 Saint \'Etienne du Rouvray, France\\}

\address[Novelli and Thibon]{Institut Gaspard Monge, Universit\'e de
Marne-la-Vall\'ee \\
5, Boulevard Descartes \\Champs-sur-Marne \\77454 Marne-la-Vall\'ee cedex 2 \\
FRANCE}
\email[Florent Hivert]{hivert@univ-rouen.fr}
\email[Jean-Christophe Novelli]{novelli@univ-mlv.fr}
\email[Jean-Yves Thibon]{jyt@univ-mlv.fr}
  
\begin{abstract}
Let $M_n$ be the $n!\times n!$ matrix indexed by permutations of $\SG_n$,
defined by $M_n(\sigma,\tau)=1$ if every descent of $\tau^{-1}$ is also a
descent of $\sigma$, and $M_n(\sigma,\tau)=0$ otherwise.
We prove the following result, conjectured by P. Dehornoy:
the characteristic polynomial $P_n(x)=|xI-M_n|$ of $M_n$ divides $P_{n+1}(x)$
in $\Z[x]$.
\vskip 0.5\baselineskip
\noindent{\bf R\'esum\'e.}
Soit $M_n$ la matrice $n!\times n!$, indexée par les \'el\'ements de $\SG_n$
et définie par $M_n(\sigma,\tau)=1$ si toute descente de $\tau^{-1}$ est aussi
une descente de $\sigma$, et $M_n(\sigma,\tau)=0$ sinon.
Nous démontrons le résultat suivant, conjecturé par P. Dehornoy :
le polynôme caractéristique $P_n(x)=|xI-M_n|$ de $M_n$ divise
$P_{n+1}(x)$ dans $\Z[x]$.
\end{abstract}

\maketitle


\section{Introduction}

On note $\SG_n$ le groupe symétrique sur $n$ éléments.
Rappelons qu'une \emph{descente} d'une permutation
$\sigma\in\SG_n$ est un entier $i$ tel que $\sigma(i)>\sigma(i+1)$.
Un \emph{recul} de $\sigma$ est une descente de son inverse
$\sigma^{-1}$. On note $\Des(\sigma)$ et $\Rec(\sigma)$ les
ensembles de descentes et de reculs de $\sigma$.

Rappelons encore que toute permutation peut s'interpréter comme une tresse
simple.
Une suite finie $(\sigma_i)_{i=1\dots l}$ de tresses simples, ou encore de
permutations, est dite \emph{normale} si et seulement si, pour tout $i<l$,
on a
\begin{equation}
\Rec(\sigma_{i+1}) \subset \Des(\sigma_{i})\,.
\end{equation}

Pour compter le nombre de suites normales de longueur $n$ et en particulier,
avoir une idée du comportement asymptotique de ce nombre quand $n\to\infty$,
Dehornoy \cite{Deh,Deh2} introduit la matrice d'adjacence du graphe dont les
chemins correspondent aux suites normales~: $M_n$ est de dimension
$n!\times n!$, avec
\begin{equation}
  M_n(\sigma, \tau ) :=
  \begin{cases}
    1 \quad\text{ si $\Rec(\tau) \subset \Des(\sigma)$,} \\
    0 \quad\text{ sinon.}
  \end{cases}
\end{equation}

\noindent
On veut montrer la conjecture suivante~\cite{Deh}~:
\begin{conjecture}\label{laconj}
Le polynôme caractéristique $P_n(x)=|xI-M_n|$ de $M_n$ divise
$P_{n+1}(x)$ dans $\Z[x]$.
\end{conjecture}

\medskip
Pour cela, nous allons interpréter la suite de matrices $(M_n)$ comme
un endomorphisme $\Phi$ de l'algèbre de Hopf $\FQSym$ des 
\og fonctions quasi-symétriques libres \fg,
et exhiber une dérivation $\delta$ qui commute avec $\Phi$.

\section{Interprétation dans les fonctions quasi-symétriques libres}

Rappelons que $\FQSym$ est une algèbre de Hopf graduée connexe \cite{NCSF6},
dont une base en degré $n$ est formée par des éléments
$\F_\sigma$, $\sigma\in\SG_n$,
qui se multiplient par \og mélange décalé \fg, c'est-à-dire, pour
$\alpha\in\SG_k$ et $\beta\in\SG_l$,
\begin{equation}
\F_\alpha \F_\beta = \sum_{\gamma\in \alpha\shuffle\beta[k]}F_\gamma
\end{equation}
où $\beta[k]$ désigne le mot dont la $i$ème lettre est $\beta_i+k$,
et $\shuffle$ le produit de mélange (ou battage) usuel.


Soit $\Phi_n$ l'endomorphisme de $\FQSym_n$ défini par
\begin{equation}
\Phi(\F_\sigma)=\sum_{\Rec(\tau)\subseteq\Des(\sigma)}\F_\tau\,.
\end{equation}
La matrice de $\Phi_n$ dans la base $\F$ est la transposée de $M_n$
et la somme directe $\Phi=\bigoplus_n\Phi_n$ est un endomorphisme
de degré $0$ de $\FQSym$.

Nous allons construire une dérivation $\partial$, surjective, de degré $-1$,
vérifiant
\begin{equation}\label{commu}
\partial\circ\Phi =\Phi\circ \partial\,.
\end{equation}

L'existence de $\partial$ vérifiant (\ref{commu}) entraîne la
conjecture~\ref{laconj}.
En effet, soit $K = \ker(\partial)$ en degr\'e $n$ et $L$ un supplémentaire
de $K$ dans $\FQSym_n$, de sorte que $\FQSym_n = K \oplus L$.
L'endomorphisme $\Phi_n$ laisse $K$ stable car
$\partial(\Phi_n(x)) = \Phi_{n-1}(\partial(x))$
et donc $\partial(x) = 0$ entra\^\i ne $\partial(\Phi_n(x))=0$.
Donc, sur la décomposition $\FQSym_n = K \oplus L$,
la matrice est triangulaire par blocs :
\begin{equation}
\begin{pmatrix}A&B\\0&C\end{pmatrix}
\end{equation}
et ainsi, le polynôme caractéristique de $\Phi_n$ est le produit de ceux de
$A$ et $C$, où $A$ est la matrice de la restriction de $\Phi_n$ à $K$
et $C$ la matrice de $\Phi_n$ sur le quotient $\FQSym_n / K$.
Par surjectivité, $\FQSym_n / K$ est isomorphe à $\FQSym_{n-1}$, 
d'où le résultat.
De plus, on peut v\'erifier que la divisibilit\'e a bien lieu dans $\Z[x]$,
car il est facile de calculer explicitement le coefficient du terme de plus
bas degr\'e de $P_n$.

\section{Construction de la dérivation équivariante}

Soit $\sigma$ un \'el\'ement de $\SG_n$ et $i$ un entier. On note $\sigma'$
le mot $(n+1)\cdot \sigma\cdot 0$.
Soit $u_i$ (resp. $v_i$) la lettre qui précède (resp. suit) la lettre $i$ dans
$\sigma'$ pour $i$ dans l'intervalle $[1,n]$.
On note alors
\begin{equation}
  \sign_i(\sigma) :=
  \begin{cases}
    +1 & \text{si $u_i<i<v_i$},\\
    -1 & \text{si $u_i>i>v_i$},\\
    0  & \text{sinon.}
  \end{cases}
\end{equation}
et $\del_i(\sigma)$ le standardisé du mot (c'est-\`a-dire la permutation ayant
les inversions aux m\^emes places que ce mot) obtenu en supprimant $i$ de
$\sigma$. Soit $\partial$ l'application linéaire définie par
\begin{equation}
  \partial_i \F_\sigma := \sign_i(\sigma) \F_{\del_i(\sigma)}
  \qquad\text{et}\qquad
  \partial := \sum_{i=1}^n \partial_i \,.
\end{equation}

\begin{lemme}
Soient $\sigma\in\SG_n$ et $\tau\in\SG_m$.
Pour tout $i\in[1,n]$, on a
\begin{equation}
\partial_i(\F_\sigma\F_\mu) = \partial_i(\F_\sigma)\ \F_\mu.
\end{equation}
Pour tout $i\in[n+1,n+m]$, on a
\begin{equation}
\partial_i(\F_\sigma\F_\mu) = \F_\sigma \partial_{i-n}(\F_\mu).
\end{equation}
\end{lemme}

{\bf D\'emonstration.}
On reprend les notations utilis\'ees pour d\'efinir $\sign$.
Soit $i\in[1,n]$.
On va consid\'erer l'ensemble $X$ des mots $w$ du m\'elange d\'ecal\'e de
$\sigma$ et $\tau$ tels que $\sign_i(w)\not=0$.
\begin{itemize}
\item Si $u_i<i<v_i$, $X$ est l'ensemble des mots du m\'elange d\'ecal\'e de
la forme
\begin{equation}
\dots u_i i B v_i \dots
\end{equation}
o\`u $B$ ne contient que des lettres de $\tau[n]$. En particulier, tous ces
mots ont pour image $1$ par $\sign_i$, et comme leurs images par $\del_i$ est
l'ensemble des permutations apparaissant dans
$\F_{\del_i(\sigma)}\F_\tau$, on conclut la d\'emonstration dans ce cas.
\item Si $u_i>i>v_i$, $X$ est l'ensemble des mots du m\'elange d\'ecal\'e de
la forme
\begin{equation}
\dots u_i B i v_i \dots
\end{equation}
o\`u $B$ ne contient que des lettres de $\tau[n]$. On conclut donc comme
dans le premier cas.
\item Si $u_i>i<v_i$, $X$ est \'egal \`a l'ensemble vide, de sorte que
\begin{equation}
\partial_i(\F_\sigma\F_\mu)=0= \partial_i(\F_\sigma)\ \F_\mu.
\end{equation}
\item Enfin, si $u_i<i>v_i$, $X$ contient les mots du m\'elange d\'ecal\'e de
la forme
\begin{equation}
\dots u_i i B v_i \dots
\qquad \text{et} \qquad
\dots u_i B i v_i \dots
\end{equation}
L'image par $\sign_i$ du premier (resp. second) ensemble est $1$ (resp. $-1$).
Comme ils ont m\^eme image par application de $\del_i$, les paires
d'\'el\'ements ayant m\^eme ensemble $B$ apportent une contribution nulle \`a 
$\partial_i(\F_\sigma\F_\mu)$, de sorte que
\begin{equation}
\partial_i(\F_\sigma\F_\mu)=0= \partial_i(\F_\sigma)\ \F_\mu.
\end{equation}
\end{itemize}
La seconde \'equation se montre de m\^eme.
\qed

\medskip
On en d\'eduit directement le r\'esultat :
\medskip

\begin{proposition}
L'application $\partial$ est une dérivation de $\FQSym$.
\end{proposition}

\medskip
Cette propri\'et\'e permet maintenant d'obtenir sans difficult\'e le
r\'esultat principal~:
\medskip

\begin{theoreme}
Les endomorphismes $\Phi$ et $\partial$ commutent.
\end{theoreme}
\medskip

{\bf D\'emonstration.}
Rappelons que si $D=\{d_1<\dots<d_k\}$ est l'ensemble des descentes d'une
permutation $\sigma\in\SG_n$, on peut recoder $D$ par la suite d'entiers
$I=(d_2-d_1,\dots,d_k-d_{k-1},n-d_k)$ appel\'ee
\emph{composition des descentes} de $\sigma$ et not\'ee $I=\CDes(\sigma)$.

Notons alors
\begin{equation}
S^I := \sum_{\Rec(\tau)\subseteq\Des(\sigma)}\F_\tau\,,
\end{equation}
o\`u $I=\CDes(\sigma)=(i_1,\dots,i_k)$.
Il est alors bien connu (\emph{cf.}~\cite{NCSF6}) que
\begin{equation}
S^I =
\F_{12..i_1}
\F_{12..i_2}
\dots
\F_{12..i_k}.
\end{equation}
Comme $\partial$ est une d\'erivation, on en d\'eduit que
\begin{equation}
\partial(\Phi(\F_\sigma)) = \sum_{I'} C_{I}^{I'} S^{I'},
\end{equation}
o\`u l'ensemble des $I'$ s'obtient en retranchant $1$ successivement \`a
chaque part $i$ de $I$ et o\`u $C_I^{I'}$ est \'egal \`a la valeur de cette
part moins deux.

On v\'erifie imm\'ediatement que $\Phi(\partial(\F_\sigma))$ est \'egal \`a
la m\^eme somme.
\qed

\medskip
Comme l'application $\partial$ est surjective, on en d\'eduit

\medskip
\begin{corollaire}[Conjecture de Dehornoy]
Le polynôme caractéristique de $M_n$ divise celui de $M_{n+1}$.
\end{corollaire}

\section{Remarques}

Les applications $\Phi$ et $\partial$ descendent à divers
quotients et sous-algèbres de $\FQSym$, en particulier
aux fonctions symétriques non-commutatives, aux fonctions
quasi-symétriques et aux fonctions symétriques ordinaires. Il est intéressant
d'observer que si on identifie $\Phi$ à un élément de
$\FQSym\otimes\FQSym^*$, il s'écrit
\begin{equation}
\Phi = \sum_{Des(\tau)\subseteq \Des(\sigma)}\F_\sigma\otimes\F_{\tau^{-1}}
=\sum_I R_I\otimes S^I=J_0(A,B)^{-1}
\end{equation}
où $J_0$ est la fonction de Bessel non-commutative introduite
dans~\cite{NT}, g\'en\'eralisant la fonction de Bessel usuelle et proposant
un rel\`evement non-commutatif des identit\'es de~\cite{Ca,CSV2}.
La commutation de $\partial$ et de $\Phi$ équivaut à une équation
fonctionnelle intéressante pour son inverse.



\begin{thebibliography}{00}
%
\bibitem{Ca} L. Carlitz,
The coefficients of the reciprocal of $J_0(z)$,
Arch. Math. 6 (1955), 121--127.
%
\bibitem{CSV2} L. Carlitz, R. Scoville, and T. Vaughan,
Enumeration of pairs of sequences by rises, falls and levels,
Manuscripta Math. 19 (1976), 211--243.
%
\bibitem{Deh} P. Dehornoy,
Combinatorics of normal sequences of braids,
J. Comb. Th. Series A, 114 (2007) 389-409; arXiv: math.CO/0511114.
%
\bibitem{Deh2} P. Dehornoy,
Still another approach to the braid ordering,
Pacific Math. J., to appear; arXiv: math.GR/0506495.
%
\bibitem{NCSF6} G. Duchamp, F. Hivert, J.-Y. Thibon,
Noncommutative symmetric functions VI: free quasi-symmetric functions and
related algebras,
Internat. J. Alg. Comput. 12 (2002), 671--717.  
%
%
%
\bibitem{NT} J.-C. Novelli, J.-Y. Thibon,
Noncommutative Bessel Symmetric Functions,
Canadian Mathematical Bulletin, to appear (also math.CO/0602043).
%
%
\end{thebibliography}
\end{document}